\documentclass[12pt]{amsart}
\usepackage{fullpage}
\usepackage{bm}
\usepackage[all]{xy}
\usepackage{tikz}
\usepackage{amsmath, amsthm, amscd, amssymb, amsfonts, latexsym}
\usepackage{mathrsfs,amssymb}
\usepackage{amsfonts}
\usepackage{amssymb}
\usepackage{amsmath}
\usepackage{amsthm}
\usepackage{bbm}
\usepackage{dsfont}
\usepackage{relsize}
\usepackage[mathscr]{euscript}
\usepackage{upgreek}
\usepackage{dsfont}
\usepackage{marvosym}

\newtheorem{theorem}{Theorem}[section]

\newtheorem{conjecture}[theorem]{Conjecture}

\newcommand{\kommentar}[1]{ }

\author{Alexandra Florea}

\address{UC Irvine, Mathematics Department, Rowland Hall, Irvine 92697, USA}
\email{floreaa@uci.edu}

\title{A Survey of Moment Bounds for $\zeta(s)$: from Heath-Brown's work to the present}

\allowdisplaybreaks

\begin{document}

\begin{abstract}
In this expository article, we review some of the ideas behind the work of Heath-Brown (D.~R. Heath-Brown, {\em Fractional moments of the Riemann Zeta Function}, J. Lond. Math. Soc., (2), {\bf 24}, no.1, (1981), 65--78) on upper and lower bounds for moments of the Riemann zeta-function, as well as the impact this work had on subsequent developments in the field. We survey recent results on the topic, which essentially recover the expected rate of growth for all moments - unconditionally for small moments and conditionally on the Riemann hypothesis for all larger moments. 

\end{abstract}

\maketitle

\section{Heath-Brown's work on fractional moments of $\zeta(s)$}
\label{hb}
This survey article is concerned with moments of the Riemann zeta-function $\zeta(s)$, focusing on the influential article \cite{hb} of Heath-Brown on fractional moments of $\zeta(s)$, and the developments that followed in the four decades since. 

Let $$\zeta(s) = \sum_{n=1}^{\infty} \frac{1}{n^s},$$ be the Riemann zeta-function,  defined for $\Re(s) >1$ and continued meromorphically to the complex plane. The moments of $\zeta(s)$  are defined as
$$I_k(T,\sigma)  = \int_1^T |\zeta(\sigma+it)|^{2k} \, dt.$$
When $\sigma=1/2$, we write $I_k(T)=I_k(T,1/2)$. 

The behavior of $I_k(T)$ is a classical topic of study since the work of Hardy and Littlewood \cite{hl}. One motivation for studying moments is that they contain information about pointwise bounds for $\zeta(1/2+it)$. Namely, showing that $I_k(T) \ll T^{1+\varepsilon}$ for all integers $k$, is equivalent to the Lindel\"of hypothesis, which states that 
\begin{equation}
\label{lh}
|\zeta(1/2+it)| \ll |t|^{\varepsilon},
\end{equation} for any $\varepsilon>0$, and $t$ large.

The Lindel\"of hypothesis is a consequence of the Riemann hypothesis (RH), one of the central problems in number theory. RH states that all the nontrivial zeros of $\zeta(s)$ lie on the vertical line $\Re(s)=1/2$. Zeros of $\zeta(s)$, and of $L$--functions more generally, play a fundamental role: many classical problems about primes can be reformulated in terms of the distribution of zeros. For example, RH is equivalent to proving the Prime Number Theorem with a strong error term, namely that the following asymptotic formula holds
$$ \pi(x) := \sum_{p \leq x} 1 = \int_{2}^x \frac{du}{\log u} + O \Big(  x^{1/2}\log x\Big),$$ as $x \to \infty$, where the sum on the left-hand side above is over prime numbers. 

A classical result due to Littlewood \cite{littlewood} shows that RH implies a strong form of the Lindel\"of hypothesis; specifically, on RH, Littlewood showed that, for all sufficiently large $t$, there exists a constant $C$ such that
$$ \Big| \zeta \Big( \frac{1}{2}+it \Big) \Big| \ll \exp \Big(C  \frac{  \log t}{\log \log t} \Big).$$
The bound above has not been improved since, except for the value of $C$; see, for example, \cite{ramachandra7, soundupperbounds}. The best known constant, $C= \frac{\log 2}{2}$, is due to Chandee and Soundararajan \cite{sound_fai}. 

The Lindel\"of hypothesis, while implied by RH, is itself a difficult open problem with significant applications to prime number distribution and to the study of zeros of $\zeta(s)$. For example, the Lindel\"of hypothesis would imply that $p_{n+1} - p_n \ll p_n^{1/2+\varepsilon}$, where $p_n$ denotes the $n^{\text{th}}$ prime number. While obtaining the bound \eqref{lh} is currently open, one can obtain weaker results, which can nevertheless have interesting applications. Unconditionally, using the functional equation of $\zeta(s)$ and the Phragm\'en-Lindel\"of principle (which provides growth bounds for analytic functions), it follows that $$\Big|\zeta \Big( \frac{1}{2} +it \Big) \Big|\ll |t|^{1/4+\varepsilon},$$  
for any $\varepsilon>0$, and large $t$. This bound is known as the convexity bound. Obtaining an exponent less than $1/4$ on the right-hand side above is a  difficult question, known as the subconvexity problem, which has attracted a lot of attention in the last century. Weyl improved the bound above to $|\zeta(1/2+it)| \ll |t|^{1/6+\varepsilon}$. Weyl’s $1/6$ exponent has since been improved, for example by Bombieri-Iwaniec \cite{bi} and Huxley \cite{huxley}. The best known subconvex bound for $\zeta(s)$ is due to Bourgain \cite{bourgain}, who showed that $|\zeta(1/2+it)| \ll  |t|^{13/84+\varepsilon}$. 

It is in the context of proving pointwise bounds for $\zeta(s)$ that Hardy and Littlewood introduced the moments $I_k(T)$ and realized that understanding the average behavior of $\zeta(s)$ provides insight into its  extreme values. Despite a long history, elucidating the structure of $I_k(T)$ has proven very difficult, and it is only in the past twenty-five years that a better understanding of the moments has emerged. 
Hardy and Littlewood \cite{hl} proved that 
$$I_1(T) \sim T \log T,$$ with a further improvement by Ingham who obtained a power-saving asymptotic formula for the second moment (which corresponds to the case $k=1$). Ingham \cite{ingham} also obtained an asymptotic for the fourth moment (corresponding to $k=2$), namely
$$I_2(T) \sim \frac{1}{2 \pi^2} T (\log T)^4.$$
The above was refined by Heath-Brown \cite{hb4}, who showed that
$$ I_2(T) = TP_4(\log T) + O(T^{7/8+\varepsilon}),$$ for some explicit polynomial $P_4$ of degree $4$. Heath-Brown \cite{hb12} also showed that $I_6(T) \ll T^{2+\varepsilon}$, which yields the subconvex Weyl-bound $|\zeta(1/2+it)| \ll t^{1/6+\varepsilon}$.

Using heuristic arguments (as explained in Section \ref{conjectures}), one expects that
\begin{equation}
\label{conj}
I_k(T) \sim C_k T(\log T)^{k^2},
\end{equation} as $T \to \infty$, for a certain constant $C_k$. At the time of Heath-Brown's 1981 paper, no plausible conjecture for the value of $C_k$ was available apart from the known values of $C_1$ and $C_2$; this changed with the work of Keating and Snaith \cite{Keating-Snaith}, as will be described in Section \ref{conjectures}. 

Since full asymptotics remain out of reach for higher moments, one naturally asks for weaker results: upper and lower bounds of the correct order of magnitude, consistent with \eqref{conj}. In other words, we are interested in proving the bounds
\begin{equation}
I_k(T) \gg T(\log T)^{k^2},
\label{lb}
\end{equation}
and 
\begin{equation}
I_k(T) \ll T (\log T)^{k^2}.
\label{ub}
\end{equation}
Before Heath-Brown \cite{hb}, Ramachandra \cite{ramachandra1, ramachandra2} had obtained the lower bound \eqref{lb} assuming that $2k$ is a positive integer. Under RH, he also showed that \eqref{lb} holds for all $k\geq 0$, and unconditionally, that
$$I_k(T) \gg T (\log T)^{k^2} (\log \log T)^{-c_k},$$
for some constant $c_k$ depending on $k$, for $k \geq 1/2$. In the other direction, Ramachandra \cite{ramachandra3} proved the sharp bound \eqref{ub} for $k=1/2$, and conditionally on RH, for $0<k<2$. 

In \cite{hb}, Heath-Brown extended these results substantially, and recovered the bounds obtained by Ramachandra, using a different, more unified approach. More precisely, he proved the following.

\begin{theorem}[Heath-Brown, \cite{hb}]
\label{main_thm}
\begin{enumerate}
\item The lower bound \eqref{lb} holds for any rational $k \geq 0$, and under RH, for any $k \geq 0$.
\item The upper bound \eqref{ub} holds for any $k =1/n$, where $n$ is a positive integer, and under RH, for any $0 \leq k \leq 2$. 
\end{enumerate}
\end{theorem}

A central ingredient in the proof of Theorem \ref{main_thm} is a convexity argument of Gabriel \cite{gabriel}, used in order to circumvent the fact that $\zeta^k(s)$ might not be regular, and hence one cannot modify contours of integration using Cauchy's theorem. 

To sketch the main ideas of proof, write $k=u/v$ where $u,v$ are positive coprime integers when $k \in \mathbb{Q}$, $u=1$ and $v$ is an integer in the case when $k=1/n$, $v=1$ and $u=k$ is any real number for the conditional lower bound in Theorem \ref{main_thm}, and finally, where $v=1$ and $1<u=k<2$ in the conditional upper bound. 

Let
$$S(s) = \sum_{n=1}^N \frac{d_k(n)}{n^s},$$ where $d_k(n)$ denotes the generalized $k^{\text{th}}$ divisor function, and where $N=T^{1/2}$ for the upper bound part of Theorem \ref{main_thm}, and $N=T^{1-\varepsilon}$ in the lower bound part, where $\varepsilon = (2-u)/4$. Let
$$g(s) = \zeta^u(s) - S^v(s),$$ 
$$w(t) = \int_T^{2T} e^{-2k(t-\tau)^2} \, d\tau,$$ and
$$J(\sigma)  = \int_{-\infty}^{\infty} |\zeta(\sigma+it)|^{2k} w(t) \, dt,$$
$$K(\sigma) = \int_{-\infty}^{\infty} |g(\sigma+it)^{2/v} w(t) \, dt,$$
$$L(\sigma) = \int_{-\infty}^{\infty} |S(\sigma+it)|^{2} w(t) \, dt.$$
Since 
$$|S^v(s)|^{2/v} = |\zeta^u(s)-g(s)|^{2/v} \ll |\zeta(s)|^{2k} + |g(s)|^{2/v},$$ it follows that
\begin{equation}
L(\sigma) \ll J(\sigma)+K(\sigma),
\label{lbound}
\end{equation} and similarly 
\begin{equation}
J(\sigma) \ll L(\sigma)+K(\sigma),
\label{jbound}
\end{equation}
\begin{equation}
\label{kbound}
K(1/2) \ll J(1/2)+L(1/2).
\end{equation}
Using the Montgomery-Vaughan mean-value theorem \cite{mv}, we have that
$$L(1/2) \asymp T (\log T)^{k^2},$$ and 
$$L(\sigma) \asymp T (\sigma-1/2)^{-k^2},$$ for $\sigma-\frac{1}{2} \gg \frac{1}{\log T}.$ 
The desired upper and lower bounds for $J(1/2)$ are obtained by combining the bounds \eqref{lbound}, \eqref{jbound}, \eqref{kbound} with delicate convexity estimates which follow from work of Gabriel \cite{gabriel}. These convexity estimates lead to upper bounds for $J(1/2)$ in terms of $J(\sigma)$ (and for $J(\sigma)$ in terms of $J(1/2)$ respectively), as well as for $K(\sigma)$ in terms of $K(1/2)$. Finally, the bounds for $I_k(T)$ follow from the corresponding bounds for $J(1/2)$, since $w(t) \ll 1$ for all $t$, and since
$$w(t) \ll \exp (-k(t^2+T^2)/18)$$ for $t \leq 0$ and $t \geq 3T$. 

We emphasize that the results proved in Heath-Brown's paper \cite{hb} not only extended the ranges in which the bounds \eqref{lb} and \eqref{ub} held and simplified existing proofs from the literature, but introduced the idea of using Gabriel's convexity estimates to address the moment questions. This idea proved to be influential, and was subsequently used by different authors in many works  to make progress on these problems. For example, concerning the moments
$$I_k(T,T+H) = \int_T^{T+H} |\zeta(1/2+it)|^{2k} \, dt,$$ 
where $T \geq H \geq 100 \log \log T \geq C_0$ (for $C_0$ a large positive constant),
Ramachandra \cite{ramachandra4, ramachandra6} showed  that
$$I_k(T,T+H) \gg H (\log H)^{k^2},$$ for any rational $k$, and that
\begin{equation}
\label{ramachandra_weaker}
I_k(T, T+H) \gg  H g \Big( \frac{\log H}{\log \log H}\Big)^{k^2},
\end{equation}
for any irrational $k$, where $g$ is a function tending to infinity with $H$. In later work \cite{hb2}, Heath-Brown, motivated by ideas in \cite{hb} and \cite{conrey_ghosh}, obtained an implied constant depending on $k$ in the bound \eqref{ub} for $0<k<2$, conditional on RH (here, the convexity arguments used in \cite{hb} were replaced by certain inequalities involving derivatives). These ideas were also used in the context of obtaining upper bounds for moments of Dirichlet $L$--functions (in the $q$-aspect) in \cite{hb_q}. Furthermore, the techniques in Heath-Brown's work have been widely applied in other contexts; for example, in studying the value-distribution of $\zeta(s)$ on the critical line \cite{jutila} (which required a uniform version in $k$ of Heath-Brown's result) or in proving upper and lower bounds in various other families of $L$--functions (in the $t$-aspect); for a non-exhaustive list, see \cite{fom1, laur1, zam, laur2, laur3, laur4}.

Since Heath-Brown's work, the area of moments has seen great developments. In the early 2000s, Keating and Snaith \cite{Keating-Snaith} made remarkable conjectures about all the moments of $\zeta(s)$ by drawing analogies between the zeta-function and characteristic polynomials of large random matrices. Progress towards these general conjectures has been achieved on multiple fronts, from proving upper bounds conditional on RH for all the moments of $\zeta(s)$, to obtaining unconditional lower bounds for all the moments and unconditional upper bounds for small $k$. We will describe the Keating-Snaith conjectures in Section \ref{conjectures}, and the subsequent progress towards them in Section \ref{progress}.

\section{Conjectures for moments of $\zeta(s)$}
\label{conjectures}

As mentioned in Section \ref{hb}, one expects the asymptotic formula in equation \eqref{conj} to hold; while the order of magnitude $T (\log T)^{k^2}$ for the $2k^{\text{th}}$ moment follows from heuristic arguments, the precise value of $C_k$ was more mysterious, and it wasn't until the breakthrough work of Keating and Snaith \cite{Keating-Snaith}, based on analogies with random matrix theory, that a clear understanding of the behavior of moments has emerged. 

When considering the shifted moments $I_k(T,\sigma)$, one sees that for $\sigma>1$, 
$$I_k(T,\sigma) \sim T \sum_{n=1}^{\infty} \frac{d_k(n)^2}{n^{2\sigma}}.$$  Focusing on the moments $I_k(T)$, note that the right-hand side in the formula above diverges for $\sigma=1/2$, but one could truncate the series to $n \leq T$, and in that case, one has that
$$ \sum_{n \leq T} \frac{d_k(n)^2}{n} \sim \frac{a_k}{\Gamma(k^2+1)} (\log T)^{k^2},$$ where 
\begin{equation}
\label{ak}
a_k= \prod_p \Big(1-\frac{1}{p} \Big)^{k^2} \Big( \sum_{j=0}^{\infty} \frac{d_k(p^j)^2}{p^j} \Big).
\end{equation}

The above computation suggests the rate of growth of the $2k^{\text{th}}$ moment given in \eqref{conj}. 
Yet another heuristic predicts the same leading order term  as in \eqref{conj}. Assuming a uniform version of Selberg's Central Limit Theorem \cite{selberg1, selberg2}, which says that $\log |\zeta(1/2+it)|$ is distributed like a Gaussian with mean $0$ and variance $\frac{1}{2} \log \log T$, and using the fact that if $X$ is a random variable with mean $\mu$ and variance $\sigma^2$, we have
$$ \mathbb{E}(e^{tX}) =  \frac{1}{ \sqrt{2 \pi } \sigma} \int_{-\infty}^{\infty} \exp \Big( tz - \frac{ (z-\mu)^2}{2 \sigma^2} \Big) \, dz = e^{t \mu + \frac{t^2 \sigma^2}{2}},$$ it would follow that
$$ I_k(T) = \int_1^T |\zeta(1/2+it)|^{2k} \, dt = \int_1^T \exp \Big( 2k \log |\zeta(1/2+it)| \Big) \, dt \asymp T (\log T)^{k^2}.$$

In light of the previous two heuristics, it is natural then to conjecture that \eqref{conj} might hold, with a value of $C_k$ yet to be determined. 
Conrey and Ghosh suggested writing $C_k = g_k a_k/\Gamma(k^2+1)$, and further proposed that $g_k$ might be an integer when $k$ is an integer. Other than the  values of $g_1$ and $g_2$ (known to be $g_1 = 1$ from \cite{hl} and $g_2=2 $ from \cite{ingham}), no plausible conjecture for the value of $g_k$ (or, equivalently, $C_k$) was put forward until the work of Conrey and Ghosh \cite{conrey_ghosh} who conjectured that $g_3=42$. Further work of Conrey and Gonek \cite{conrey_gonek} recovered the conjecture for $g_3$ and predicted that $g_4 = 24024$. Extending the Conrey and Gonek heuristic arguments to higher moments did not yield conjectures for $k>4$. The breakthrough work of Keating and Snaith \cite{Keating-Snaith}, based on analogies with random matrix theory, produced the following very general conjecture.
\begin{conjecture}
\label{conjecture_ks}
Let $k$ be a positive real number. Then
$$\lim_{T \to \infty} \frac{I_k(T)}{T (\log T)^{k^2}} = \frac{g_k a_k}{\Gamma(k^2+1)},$$ with $a_k$ given in \eqref{ak}, and
$$g_k = \Gamma(k^2+1) \frac{G^2(1+k)}{G(1+2k)},$$
where $G$ denotes the Barnes $G$-function.
\end{conjecture}
We note that for integer $k$, the expression above simplifies to
$$g_k = (k^2)! \prod_{j=0}^{k-1} \frac{j!}{(j+k)!}.$$

The main idea behind the Keating and Snaith conjecture was to model $\zeta(s)$ by the characteristic polynomial of random unitary matrices. Denote by $1/2+it_n$ the non-trivial zeros of $\zeta(s)$. Under RH, it follows that $t_n \in \mathbb{R}$. The motivation behind the Keating-Snaith conjecture stands in Montgomery's conjecture \cite{montgomery} that the two-point correlations between the heights $t_n$ of the non-trivial zeros of $\zeta(s)$, as $n \to \infty$, are the same as those between the eigenvalues of random complex hermitian matrices as the size of the matrix goes to infinity. These matrices form the Gaussian Unitary Ensemble (GUE) in random matrix theory, and their eigenvalues have the same local statistics as the unitary matrices in the limit as the matrix size tends to infinity. Montgomery proved the conjecture in certain ranges, and there is also strong numerical evidence supporting it \cite{odlyzko}. The available evidence thus suggests that, as $n \to \infty$, the local statistics of the scaled zeros of $\zeta(s)$, denoted by $w_n = \frac{t_n}{2\pi} \log \frac{t_n}{2\pi}$ (scaled to have mean spacing one) are the same as those of normalized eigenphases $\phi_n = \theta_n \frac{N}{2 \pi}$, where $e^{i\theta_n}$ denote
the eigenvalues of $N \times N$ unitary matrices, in the limit as $N \to \infty$. This suggests that the properties of $\zeta(s)$ might be modelled by the properties of 
$$Z(U) = \det(I-U),$$ as $U$ ranges over $N \times N$ unitary matrices. Identifying the mean density of the eigenangles $\theta_n$, which is $N/2\pi$, with the mean density of the zeros of $\zeta(s)$ at height $T$, which is $\frac{1}{2 \pi } \log \frac{T}{2\pi}$, suggests setting
$N = \log \frac{T}{2\pi}.$
In order to understand the behavior of $I_k(T)$, Keating and Snaith looked at the analogous moments on the random matrix side; namely they considered
$$ \int_{U(N)} |Z(U)|^{2k} \, dU,$$ which can be evaluated using the Weyl integration formula and the Selberg integral \cite{mehta}. These give that 
$$\int_{U(N)} |Z(U)|^{2k} \, dU  \sim g_k \frac{N^{k^2}}{\Gamma(k^2+1)},$$ with $g_k$ as in Conjecture \ref{conjecture_ks}. For $k$ a positive integer, in fact Keating and Snaith proved that the integral above is equal to a polynomial in $N$ of degree $k^2$. 

While this computation motivates the factor $g_k$ in Conjecture \ref{conjecture_ks}, the factor $a_k$, which contains information about the primes, had to be put in ``by hand''. A different approach, combining the information about the primes and about the zeros, was described by Gonek, Hughes, Keating \cite{hybrid}, recovering Conjecture \ref{conjecture_ks}. Yet another model was developed by Conrey, Farmer, Keating, Rubinstein and Snaith \cite{cfkrs}, who conjectured that, for integer $k$,
$$I_k(T) = \int_0^T P_k(\log t/2\pi) \, dt + O(T^{1-\delta}),$$ for some $\delta>0$, where $P_k$ is an explicit polynomial of degree $k^2$, whose leading coefficient is equal to $g_k a_k/(k^2)!$. 

The above conjectures generalize to many other families of $L$--functions. Based on function field analogues, Katz and Sarnak \cite{KS} conjectured that the distribution of zeros near $1/2$ in families of $L$--functions matches that of eigenvalues near $1$ of random matrices in $U(N), USp(2N), SO(2N)$ or $SO(2N+1)$. The Keating-Snaith conjectures can then be formulated for various other families of $L$--functions \cite{ks2}, by considering the random matrix moments in the appropriate compact group suggested by the work of Katz and Sarnak. A different approach based on multiple Dirichlet series was developed in \cite{dgh} for the family of quadratic Dirichlet $L$--functions. Ideas which involve choosing long Dirichlet approximations to $\zeta(s)$ led to the same conjectures for moments of $\zeta(s)$ in later work of Conrey and Keating \cite{ck1, ck2, ck3, ck4, ck5}.

\section{Progress towards the Keating-Snaith conjectures}
\label{progress} 

While the Keating–Snaith Conjecture \ref{conjecture_ks} provides a compelling model for the moments of $\zeta(s)$, proving such asymptotics for all $k$ remains beyond reach. Instead, progress has focused on obtaining upper and lower bounds for $I_k(T)$ that are consistent with the conjecture. Progress follows three directions: (1) obtaining upper bounds on RH for high moments; (2) obtaining unconditional upper bounds for ``small'' moments; (3) obtaining unconditional lower bounds for all the moments. The ideas behind these results are sieve theoretic inspired; while the proofs all use ingenious different arguments, they have a unifying feel. 
\subsection{Conditional upper bounds}
Soundararajan \cite{soundupperbounds} obtained almost sharp upper bounds for all the moments $I_k(T)$, conditional on RH. Namely, he showed that for any $\varepsilon>0$ and positive real $k$, we have
\begin{equation}
\label{sound_bound}
I_k(T) \ll T(\log T)^{k^2+\varepsilon}.
\end{equation}
 Radziwi\l\l \,\cite{436} strengthened the bound above for $0<k<\frac{24}{11}$ under RH, removing the $\varepsilon$ on the power of $\log T$. This extended the range of $k$ due to Heath-Brown for which the sharp bound \eqref{ub} holds (on RH).

The almost sharp upper bound \eqref{sound_bound} was subsequently refined by Harper \cite{harper} who obtained an upper bound of the expected order of magnitude (as in equation \eqref{ub}). One of the key ideas in \cite{soundupperbounds} and \cite{harper} is bounding $\log |\zeta(1/2+it)|$ by a sum over the primes, under RH. For example, Soundararajan, going back to ideas of Selberg \cite{selberg2}, showed that, under RH, for $T \leq t \leq 2T$ and for $2 \leq x \leq T^2$, 
\begin{equation}
\label{log_ineq}
 \log | \zeta (\tfrac{1}{2}+it)| \leq \Re \sum_{n \leq x} \frac{\Lambda(n)  \log(x/n)}{n^{\frac{1}{2}+\frac{1}{\log x}+it} \log x} + \frac{ \log T}{\log x} +O \Big( \frac{1}{\log x} \Big).
 \end{equation}

If $x$ is sufficiently small, then one can compute moments of the above sum over the primes; however, the $\log T/\log x$ term might be large. Soundararajan split the sum over the primes into two ranges. On one hand, for the small primes, one can compute more moments, which allows for a better understanding. On the other hand, for the larger primes, although one can compute fewer moments, the fact that the variance of the sum over the primes is not too large provides control on the size of $\log|\zeta(1/2+it)|$. Soundararajan exploited these ideas to investigate the frequency of large values of the Dirichlet polynomial in \eqref{log_ineq} by using Markov-inequality type arguments applied to their high moments, establishing the bound \eqref{sound_bound}. 

Harper further developed these ideas in the beautiful paper \cite{harper}, but worked directly with moments rather than studying large deviations as in Soundararajan's work \cite{soundupperbounds}. Starting similarly as in \cite{soundupperbounds}, Harper bounded $\log |\zeta(1/2+it)|$ by a long Dirichlet polynomial (over the primes), and split the polynomial into multiple pieces (as opposed to two pieces as in \cite{soundupperbounds}), raising each piece to a different power. The larger primes are expected to contribute increasingly little, and they can be raised to smaller powers. On one hand, when the Dirichlet polynomial pieces are suitably small, their exponentials are approximated by suitably truncated Taylor expansions. On the other hand, one exploits the fact that if $t$ is a point for which a certain Dirichlet piece is too large, the set of such $t$ has small measure, and hence does not affect the upper bound.

To be more specific, Harper split the Dirichlet polynomial on the right hand side of \eqref{log_ineq} into primes in intervals $(T^{\beta_0}, T^{\beta_1}],(T^{\beta_1}, T^{\beta_2}],\ldots, (T^{\beta_{K-1}}, T^{\beta_K}]$, where $$\beta_0 = 0,  \beta_j = \frac{20^{j-1}}{(\log\log T)^2},$$ for $1 \leq j \leq K$, where $K$ is such that $\beta_K$ is a small constant (say $\beta_K \sim e^{-1000k}$). Let
$$\mathcal{P}_j(t) =  \sum_{T^{\beta_{j-1}} < p \leq T^{\beta_j}} \frac{1}{p^{\frac{1}{2}+it}}.$$
We note that, for the sake of exposition, the actual expression of $\mathcal{P}_j(t)$ is written in a simplified form here, and from now on we disregard the contribution of the square of the primes in the Dirichlet polynomial in \eqref{log_ineq}.

One of the key ideas in the proof is the flexibility in the choice of $x$ in \eqref{log_ineq}. Let $\mathcal{T}_0$ denote the set of exceptional points $t$ for which $|\Re \mathcal{P}_1(t)| \geq \beta_1^{-3/4}$. Using Markov's inequality, it follows that
$$ \text{meas}(\mathcal{T}_0) \beta_1^{-3\ell/2} \leq \int_0^T | \mathcal{P}_1(t)|^{2\ell} \ll T \ell! (\log \log T)^{\ell},$$ for $\ell$ suitably chosen so that one can compute the moments of $\mathcal{P}_1(t)$ above. Taking $\ell = [(\log \log T)^2/5]$, it then follows that the measure of $\mathcal{T}_0$ is negligible; namely
$$  \text{meas}(\mathcal{T}_0) \ll T e^{-(\log \log T)^2/10},$$
and hence 
\begin{align*}
\int_{\mathcal{T}_0} |\zeta(1/2+it)|^{2k} \, dt \leq I_{2k}(T)^{1/2} \text{meas}(\mathcal{T}_0)^{1/2} \ll T,
\end{align*}
where one can use a priori upper bounds for $I_{2k}(T)$ (for example, as in \cite{soundupperbounds}). Now assume that for each $r \leq K$, $|\Re \mathcal{P}_j(t)| < \beta_j^{-3/4}$ for $j \leq r$, but that $|\Re \mathcal{P}_{r+1}(t)| \geq \beta_{r+1}^{-3/4}$. Call $\mathcal{T}_r$ the set of such $t$. Then using \eqref{log_ineq} with $x=T^{\beta_r}$, it follows that
\begin{align}
\label{tr}
\int_{\mathcal{T}_r} |\zeta(1/2+it)|^{2k} \, dt \ll e^{2k/\beta_r} \int_0^T (|\Re \mathcal{P}_{r+1}(t)| \beta_{r+1}^{3/4})^{\ell_{r+1}} \prod_{j=1}^r E_{\ell_j}(2k \Re \mathcal{P}_j(t)) \, dt,
\end{align}
where $$E_{\ell}(t) = \sum_{s \leq \ell} \frac{t^s}{s!}$$ denotes a truncated Taylor series. We note that in the above, one uses the fact that $\exp(2k \Re \mathcal{P}_j(t))$ can be approximated by a suitably truncated Taylor series. The parameters $\ell_j$ are such that $\ell_j \asymp \beta_j^{-1}$, and they are chosen so that the length of the Dirichlet polynomial on the right hand side is small enough, and hence one can use mean-value theorems to approximate it. Although the term $\exp(2k/\beta_r)$ in \eqref{tr} can be big, one still wins because $\mathcal{P}_{r+1}(t)$ is rarely large, and hence one obtains a suitable bound for \eqref{tr}. 

\subsection{Unconditional upper bounds}
Following Heath-Brown's proof of the unconditional bound \eqref{ub} for $k=1/n$, where $n$ is a positive integer, unconditional sharp upper bounds were obtained by Bettin, Chandee and Radziwi\l\l \,\cite{bcr} for $k=1+1/n$, for $n$ a positive integer. 

More recently, Heap, Radziwi\l\l, and Soundararajan \cite{hrs} proved \eqref{ub} for all real $0 \leq k \leq 2$. Their proof illustrated an upper bound principle, formulated by Radziwi\l\l \,and Soundararajan \cite{rs}, which roughly says that if one can compute a $k^{\text{th}}$ moment in a family of $L$--functions with a little room to spare (i.e., with a power savings error term, or equivalently, a $k^{\text{th}}$ moment multiplied by a short Dirichlet polynomial), then one can obtain sharp upper bounds for all the moments less than $k$. 
Proving \eqref{ub} unconditionally for $0 \leq k \leq 2$ uses an iterative scheme similar in flavor to the one described above coming from Harper's work \cite{harper} (and bearing similarities to arguments developed in \cite{rs}), together with knowledge of the fourth moment of zeta twisted by a short Dirichlet polynomial (as in the work of Hughes and Young \cite{hughes_young} or Bettin, Bui, Li and Radziwi\l\l \, \cite{bblr}).

\subsection{Unconditional lower bounds}
Following the unconditional lower bound \eqref{lb} for all rational $k$ proven in \cite{hb}, and Ramachandra's weaker bound $I_k(T) \gg T (\log T)^{k^2} (\log \log T)^{-k^2}$ (see equation \eqref{ramachandra_weaker}), there has been significant progress on obtaining sharp lower bounds.

Rudnick and Soundararajan \cite{rudnick_sound} developed an  approach that yielded tight lower bounds for rational moments in various families of $L$--functions. Their method was subsequently refined by 
Radziwi\l\l \, and Soundararajan \cite{rs_lb} to produce sharp lower bounds for all $k \geq 1$.   

Finally, work of Heap and Soundararajan \cite{heap_sound} obtained tight lower bounds for all $0<k \leq 1$, settling the lower bound problem. As in \cite{rs_lb}, the approach can be extended to other families of $L$--functions. The lower bound principle illustrated in \cite{heap_sound} and \cite{rs_lb} can roughly be formulated as follows: if one can compute the mean value multiplied by a short Dirichlet polynomial, then one can get sharp lower bound for all the moments bigger than the first. Additionally, if one can compute the second moment in the family multiplied by a short Dirichlet polynomial, then one can recover tight lower bounds for moments below the first as well. In different work, by obtaining lower bounds for the large deviations in Selberg's Central Limit Theorem, Arguin and Bailey \cite{ab} recovered unconditional sharp lower bounds for all the positive moments of $\zeta(s)$.  

This circle of ideas has been fruitfully used in other contexts as well, to obtain bounds for shifted moments of $\zeta(s)$ \cite{chandee_shifted, curran, curran2, ng}, to study negative moments of $\zeta(s)$ \cite{bf_negative}, as well as in diverse applications to, for example, sign changes of Fourier coefficients of half-integral weight modular forms \cite{lr}, non-vanishing of $L$--functions \cite{dfl2}, or weighted central limit theorems for central values of $L$--functions \cite{buietall}.

\subsection*{Acknowledgements} The author would like to thank the referee for helpful comments and suggestions. This work was supported by the National Science Foundation (DMS-2101769 and CAREER grant DMS-2339274). 

 \bibliographystyle{amsalpha}

\bibliography{bibliography}
\end{document}